\documentclass[12pt]{article}
\usepackage{amssymb}

\newtheorem{theorem}{Theorem}[section]

\newtheorem{definition}[theorem]{Definition}

\newtheorem{proposition}[theorem]{Proposition}

\newenvironment{proof}[1][Proof]{\textbf{#1.} }{\ \rule{0.5em}{0.5em}}

\begin{document}

\title{Metric nonlinear connections}
\author{\textsc{Ioan Bucataru} \\
{\small Faculty of Mathematics, ``Al.I.Cuza'' University,} \\
{\small Ia\c{s}i, 700506, Romania, email: bucataru@uaic.ro}}
\date{}
\maketitle

\noindent\small{\textbf{Abstract} For a system of second order
differential equations we determine a nonlinear connection that is
compatible with a given generalized Lagrange metric. Using this
nonlinear connection, we can find the whole family of metric
nonlinear connections that can be associated with a system of SODE
and a generalized Lagrange structure. For the particular case when
the system of SODE and the metric structure are Lagrangian, we
prove that the canonic nonlinear connection of the Lagrange space
is the only nonlinear connection which is metric and compatible
with the symplectic structure of the Lagrange space. The metric
tensor of the Lagrange space determines the symmetric part of the
nonlinear connection, while the symplectic structure of the
Lagrange space determines the skew-symmetric part of the nonlinear
connection.}

\vspace*{3mm}

\noindent\textbf{2000 MSC: 53C05, 53C60, 34A26}

\noindent\textbf{Keywords}: semispray, nonlinear connection,
metric structure, symplectic structure, Lagrange space

\section*{Introduction}

Nonlinear connections and metric structures are important tools
for the differential geometry of the tangent bundle. Using the
dynamical covariant derivative one can associate to a nonlinear
connection, we introduce a compatibility condition between a
nonlinear connection and a generalized Lagrange metric. This
compatibility condition is a natural generalization of the well
known metric compatibility of a linear connection in a Riemannian
space, \cite{[docarmo]}.

For the differential geometry of a system of SODE one can
associate a nonlinear connection and the corresponding dynamical
covariant derivative. Such nonlinear connections were introduced
by M. Crampin \cite{[crampin1]} and J. Grifone \cite{[grifone]}. A
metric geometry of a system of SODE requires a nonlinear
connection which is compatible with a given generalized Lagrange
metric. If $S$ is an SODE and $g$ a generalized Lagrange metric,
we determine a metric nonlinear connection that corresponds to the
pair $(S,g)$. Using this nonlinear connection, we determine the
whole family of metric nonlinear connections that correspond to
the pair $(S,g)$.

For the particular case of a Lagrange space, the metric
compatibility and the compatibility with the symplectic structure
of the Lagrange space uniquely determine the nonlinear connection
one can associate with Euler-Lagrange equations. The compatibility
of the nonlinear connection with the symplectic structure of the
Lagrange space is equivalent with the existence of an almost
Hermitian structure on $TM$. The compatibility with the Lagrange
metric determines the symmetric part of the nonlinear connection.
The compatibility with the symplectic structure of the Lagrange
space determines the skew-symmetric part of the nonlinear
connection.

The metric compatibility of a semispray and the associated
nonlinear connection with a generalized Lagrange metric has been
studied by M. Crampin et al. \cite{[crampin2]}, O. Krupkova
\cite{[krupkova]}, W. Sarlet \cite{[sarlet]}, J. Szilasi and Z.
Muzsnay \cite{[szilasi2]} and it is known as one of the Helmholtz
condition for the inverse problem of Lagrangian Mechanics. In the
above mentioned papers, the problem that is studied is as follows:
for a given semispray and the associated nonlinear connection find
if it exists a Lagrange metric with respect to which the nonlinear
connection is compatible. In our work a system of SODE and a
generalized Lagrange metric are given a priori and we associate to
these structures a metric nonlinear connection. This nonlinear
connection, in general, is different from the nonlinear connection
one usually associates to a semispray. However, the two nonlinear
connections coincide if the metric structure is Lagrangian. A
geometric theory of the pair $(S,g)$ has been proposed also by B.
Lackey in \cite{[lackey]}, using Cartan's method of equivalence. A
different approach for studying metrizable nonlinear connection
has been proposed by M. Anastasiei \cite{[anastasiei]}. However,
this approach coincides with ours only for the particular case of
a Finsler space.

\section{Geometric structures on tangent bundles.}

In this section we introduce the geometric structures we deal with
in this paper: semisprays, nonlinear connections and metric
structures. These structures are defined on the total space of a
tangent bundle.

We start by considering $M$ a real, $n$-dimensional manifold of
$C^{\infty}$-class and denote by $(TM,\pi,M)$ its tangent bundle.
If $(U, \phi=(x^i))$ is a local chart at $p\in M$ from a fixed
atlas of $C^{\infty}$-class, then we denote by $(\pi^{-1}(U),
\Phi=(x^i, y^i))$ the induced local chart at $u\in \pi^{-1}(p)
\subset TM$. The linear map $\pi_{*,u}:T_uTM \rightarrow
T_{\pi(u)}M$ induced by the canonical submersion $\pi$ is an
epimorphism of linear spaces for each $u\in TM$. Therefore, its
kernel determines a regular, $n$-dimensional, integrable
distribution $V:u\in TM \mapsto V_uTM:=\mathrm{Ker} \pi_{*,u}
\subset T_uTM$, which is called the \textit{vertical
distribution}. For every $u \in TM$, $\{{\partial}/{\partial
y^i}|_u\}$ is a basis of $V_uTM$, where $\{{\partial}/{\partial
x^i}|_u, {\partial}/{\partial y^i}|_u\}$ is the natural basis of
$T_uTM$ induced by a local chart. Denote by $\mathcal{F}(TM)$ the
ring of real-valued functions over $TM$ and by $\mathcal{X}(TM)$
the $\mathcal{F}(TM)$-module of vector fields on $TM$. We shall
consider also $\mathcal{X}^v(TM)$ the $\mathcal{F}(TM)$-module of
vertical vector fields on $TM$. An important vertical vector field
is $\mathbb{C}=y^i({\partial}/{\partial y^i})$, which is called
the \textit{Liouville vector field}.

The mapping $J: \mathcal{X}(TM) \rightarrow \mathcal{X}(TM)$ given
by $J=({\partial}/{\partial y^i}) \otimes dx^i $ is called the
\textit{tangent structure} and it has the following properties:
Ker $J$ = Im $J$ = $\mathcal{X}^v(TM)$; rank $J=n$ and $J^2=0$.

A vector field $S\in\chi(TM)$ is called a \textit{semispray}, or a
second order vector field, if $JS=\mathbb{C}$. In local
coordinates a semispray can be represented as follows:
\begin{equation}
S=y^i\frac{\partial}{\partial x^i}-2G^i(x,y)
\frac{\partial}{\partial y^i}. \label{smispray}
\end{equation}
We refer to the functions $G^i(x,y)$ as to the local coefficients
of the semispray $S$. Integral curves of a semispray $S$ are
solutions of the following system of SODE:
\begin{equation}
\frac{d^2x^i}{dt^2}+2G^i\left(x,\frac{dx}{dt}\right)=0.
\label{sode}
\end{equation}
A \textit{nonlinear connection} $N$ on $TM$ is an $n$-dimensional
distribution $N: u\in TM \mapsto N_uTM\subset T_uTM$ that is
supplementary to the vertical distribution. This means that for
every $u\in TM$ we have the direct decomposition:
\begin{equation}
T_uTM=N_uTM \oplus V_uTM. \label{thv} \end{equation} The
distribution induced by a nonlinear connection is called the
\textit{horizontal distribution}. We denote by $h$ and $v$ the
horizontal and the vertical projectors that correspond to the
above decomposition and by $\mathcal{X}^h(TM)$ the
$\mathcal{F}(TM)$-module of horizontal vector fields on $TM$. For
every $u=(x,y)\in TM$ we denote by ${\delta}/{\delta
x^i}|_u=h({\partial}/{\partial x^i}|_u).$ Then $\{{\delta}/{\delta
x^i}|_u, {\partial}/{\partial y^i}|_u\}$ is a basis of $T_uTM$
adapted to the decomposition (\ref{thv}). We call it the
\textit{Berwald basis} of the nonlinear connection. With respect
to the natural basis $\{{\partial}/{\partial x^i}|_u,
{\partial}/{\partial y^i}|_u\}$ of $T_uTM$, the horizontal
components of the Berwald basis have the expression:
\begin{equation}
\left.\frac{\delta}{\delta x^i}\right|_u =
\left.\frac{\partial}{\partial x^i}\right|_u -
N^j_i(u)\left.\frac{\partial}{\partial y^j}\right|_u, \ u \in TM.
\label{coefficients}\end{equation} The functions $N^i_j(x,y)$,
defined on domains of induced local charts, are called the
\textit{local coefficients} of the nonlinear connection. The dual
basis of the Berwald basis is $\{dx^i, \delta y^i=dy^i +
N^i_jdx^j\}$.

It has been shown by M. Crampin \cite{[crampin1]} and J. Grifone
\cite{[grifone]} that every semispray determines a nonlinear
connection. Local coefficients of the induced nonlinear connection
are $N^i_j={\partial G^i}/{\partial y^j}$.

A \textit{generalized Lagrange metric}, or a GL-metric for short,
is a (2,0)-type symmetric d-tensor field $g=g_{ij}(x,y)dx^i\otimes
dx^j$ of rank $n$ on $TM$. Throughout this paper by a d-tensor
field we mean a tensor field on $TM$, whose components, under a
change of coordinates on $TM$, behave like the components of a
tensor field on the base manifold $M$. One can use the generalized
Lagrange metric to define a metric structure on the vertical
subbundle $VTM$, that is we can consider $g^v=g_{ij}\delta
y^i\otimes \delta y^j$. Then, $G=g+g^v$ is a metric structure on
$TM$ with respect to which the horizontal and the vertical
distributions are orthogonal to each other.

The geometry of the pair $(M, g_{ij}(x,y))$ is called the geometry
of a generalized Lagrange space. This geometry has been studied by
R. Miron in \cite{[miron1]}, \cite{[mirona1]}. However, in this
work no compatibility condition between the generalized Lagrange
metric and a nonlinear connection is required.

\section{Metric nonlinear connections}

Nonlinear connections, semisprays and metric structures are
important tools in the geometry of tangent bundles. There are
situations, as in the geometry of generalized Lagrange spaces,
\cite{[miron1]}, where these structures are considered, but no
condition of compatibility is required for them. Using the
covariant derivative one can associate to a semispray $S$ and a
nonlinear connection $N$, we introduce a compatibility condition
between the pair $(S,N)$ and a generalized Lagrange metric $g$.
This compatibility is a natural generalization of the well known
metric compatibility of a linear connection in a Riemannian space,
\cite{[docarmo]}. As the metric compatibility is not enough to
determine the Levi-Civita connection of a Riemannian space,
similarly the metric compatibility does not uniquely determine a
nonlinear connection. A whole family of metric nonlinear
connections is determined when a generalized Lagrange metric and a
semispray are fixed. The problem of compatibility between a system
of second order differential equations and a metric structure has
been studied by numerous authors, \cite{[lackey]},
\cite{[crampin2]} and it is known as one of the Helmholtz
conditions from the inverse problem of Lagrangian Mechanics,
\cite{[crampin2]}, \cite{[krupkova]}, \cite{[sarlet]},
\cite{[szilasi2]}. In this section we approach the Helmholtz
condition in a different way: for a given semispray and a
generalized Lagrange metric we determine the whole family of
nonlinear connections that are compatible with the metric tensor.

Let $N$ be a nonlinear connection with local coefficients
$N^i_j(x,y)$ and let $S$ be a semispray. We determine the whole
family of nonlinear connections one can associate to the semispray
$S$ and that are compatible with a generalized Lagrange metric
$g$. The dynamical covariant derivative that corresponds to the
pair $(S,N)$ is defined by $\nabla: \chi^v(TM) \longrightarrow
\chi^v(TM)$ through:
\begin{equation}
\nabla \left(X^i\frac{\partial}{\partial y^i}\right)=\left(S(X^i)+
X^jN^i_j\right) \frac{\partial}{\partial y^i}. \label{covariant}
\end{equation}
In terms of the natural basis of the vertical distribution we have
\begin{equation}
\nabla \left(\frac{\partial}{\partial y^i}\right)=N^j_i
\frac{\partial}{\partial y^j}. \label{ncovariant}
\end{equation}
Hence, $N^i_j$ are also local coefficients of the dynamical
covariant derivative. Dynamical covariant derivative $\nabla$
corresponds to the covariant derivative $D$ in \cite{[crampin2]}
or $\mathcal{D}_S$ in \cite{[krupkova]}, along the integral curves
of the semispray $S$. Dynamical covariant derivative $\nabla$ has
the following properties:
\begin{itemize}
\item[1)] $\nabla(X+Y)=\nabla X + \nabla Y, \forall X, Y \in
\chi^v(TM)$, \vspace*{1mm} \item[2)] $\nabla(fX)=S(f) X + f\nabla
X, \forall X \in \chi^v(TM), \forall f\in \mathcal{F}(TM)$.
\end{itemize}
It is easy to extend the action of $\nabla$ to the algebra of
d-tensor fields by requiring for $\nabla$ to preserve the tensor
product. For a GL-metric $g$, which is a (2,0)-type d-tensor
field, its dynamical covariant derivative is given by
\begin{equation}
(\nabla g)(X,Y)=S(g(X,Y)) - g(\nabla X,Y) - g(X,\nabla Y), \forall
X,Y. \label{covariantg}
\end{equation}
In local coordinates, we have:
\begin{equation}
g_{ij|}:=(\nabla g)\left(\frac{\partial}{\partial y^i},
\frac{\partial}{\partial y^j}\right)=S(g_{ij}) - g_{im}N^m_j -
g_{mj}N^m_i. \label{covariantgij}
\end{equation}
\begin{definition} \label{metricnon} Let $S$ be a semispray, $N$ a nonlinear connection
and $\nabla$ the associated covariant derivative. The nonlinear
connection $N$ is metric or compatible with the metric tensor $g$
if $\nabla g=0$, which is equivalent to:
$$
S(g(X,Y)) = g(\nabla X,Y) + g(X,\nabla Y), \forall X,Y.
$$
\end{definition}
For a semispray $S$ with local coefficients $G^i$ one can
associate a nonlinear connection with local coefficients
$N^i_j={\partial G^i}/{\partial y^j}$. In general this nonlinear
connection is not metric with respect to $g$. For the fixed
semispray $S$, we determine first a nonlinear connection that is
metric with respect to $g$ and then we determine the whole family
of nonlinear connections with this property.

Let us consider the following Obata operators one can associate to
a GL-metric $g_{ij}$, \cite{[mirona1]}:
\begin{equation}
O^{ij}_{kl}=\frac{1}{2}(\delta^i_k\delta^j_l - g^{ij}g_{kl})
\textrm{ and } O^{\ast ij}_{kl}=\frac{1}{2}(\delta^i_k\delta^j_l +
g^{ij}g_{kl}). \label{GLobata}
\end{equation}
\begin{theorem} \label{metricn} Let $S$ be a semispray with local coefficients
$G^i$. There is a metric nonlinear connection $N^c$, whose
coefficients $N^{ci}_j$ are given by:
\begin{equation}
N^{ci}_j=\frac{1}{2}g^{ik}S(g_{kj})+O^{ik}_{sj}\frac{\partial
G^s}{\partial y^k}. \label{metricnon1}
\end{equation}
\end{theorem}
\begin{proof} One can write coefficients $N^{ci}_j$ from
expression (\ref{metricnon1}) into the following equivalent form
\begin{equation}
N^{ci}_j=\frac{1}{2}g^{ik}g_{kj|} + \frac{\partial G^i}{\partial
y^j}. \label{metricnon2}
\end{equation}
In expression (\ref{metricnon2}) the covariant derivative
$g_{kj|}$ is taken with respect to the pair $(G^i, {\partial
G^i}/{\partial y^j})$. Since ${\partial G^i}/{\partial y^j}$ are
local coefficients of a nonlinear connection and $g^{ik}g_{kj|}$
are components of a d-tensor field of (1,1)-type we have that
$N^{ci}_j$ are also the local coefficients of a nonlinear
connection. Consider now the covariant derivative $\nabla$ one can
associate to the pair $(G^i, N^{ci}_j)$. It is a straightforward
calculation to check that
$$
S(g_{ij})-g_{im}N^{cm}_j-g_{mj}N^{cm}_i=0,
$$
which means that $\nabla g=0$ and hence, the nonlinear connection
$N^c$ is metric.
\end{proof}

Local coefficients of the metric nonlinear connection given by
expression (\ref{metricnon1}) can be written as follows:
\begin{equation}
N^{ci}_j=\frac{1}{2}g^{ik}S(g_{kj}) +
\frac{1}{2}\left(\frac{\partial G^i}{\partial y^j} - g^{ik}g_{mj}
\frac{\partial G^m}{\partial y^k}\right), \label{lackey}
\end{equation}
which coincides with the nonlinear connection determined by B.
Lackey in \cite{[lackey]}.

\begin{proposition} Let $S$ be a semispray with local coefficients $G^i$, $N$ the
associated nonlinear connection with local coefficients
$N^i_j={\partial G^i}/{\partial y^j}$, and $N^c$ the metric
nonlinear connection given by expression (\ref{metricnon1}). The
nonlinear connection $N$ is metric if and only if $N=N^c$.
\end{proposition}

The metric compatibility of the nonlinear connection
$N^i_j={\partial G^i}/{\partial y^j}$ reads as follows:
\begin{equation}
S(g_{ij})-g_{im}\frac{\partial G^m}{\partial y^j} -
g_{mj}\frac{\partial G^m}{\partial y^i}=0, \label{helmholtz}
\end{equation}
which is one of the Helmholtz conditions for the inverse problem
in Lagrangian Mechanics, \cite{[sarlet]}.

Now, we can determine the whole family of metric nonlinear
connections one can associate to a semispray.
\begin{theorem} \label{family}
Consider $S$ a semispray with local coefficients $G^i$ and $N^c$
the metric nonlinear connection with local coefficients given by
expression (\ref{metricnon1}). The family of all nonlinear
connections that are metric with respect to the GL-metric $g_{ij}$
is given by:
\begin{equation}
N^i_j=N^{ci}_j + O^{ki}_{jm} X^m_k, \label{GLmetricnon}
\end{equation}
where $X^m_k$ is an arbitrary (1,1)-type d-tensor field.
\end{theorem}
\begin{proof} The condition that both nonlinear connections
$N^{ci}_j$ and $N^i_j$ are metric with respect to the metric
tensor $g$, can be written as $S(g_{ij})=g_{mj}N^{cm}_i +
g_{im}N^{cm}_j$ and $S(g_{ij})=g_{mj}N^{m}_i + g_{im}N^{m}_j$. If
we subtract these two equations we obtain $O^{\ast is}_{jm}
(N^{m}_i-N^{cm}_i)=0$. Using the fact that Obata operators
(\ref{GLobata}) are projectors, which implies $O^{ij}_{kl}O^{\ast
km}_{pj}=0$, we obtain that the solution of this tensorial
equation is given by expression (\ref{GLmetricnon}).
\end{proof}

It is possible to define a dynamical covariant derivative $\nabla$
given by expression (\ref{covariant}) by considering a nonlinear
connection $N$, only, without considering an arbitrary semispray
$S$. In such a case the semispray $S$ is the horizontal vector
field $S=y^i({\delta}/{\delta x^i})$ with local coefficients
$2G^i(x,y)=N^i_j(x,y)y^j$. All results obtained in this section
can be reformulated within the new particular framework. However,
this does not allow us to determine a canonic metric nonlinear
connection for a generalized Lagrange space.

There are classes of generalized Lagrange spaces that posses
canonic nonlinear connections. However, these nonlinear
connections are not compatible with the generalized Lagrange
metric. Such spaces are called regular generalized Lagrange spaces
and they were introduced by R. Miron in \cite{[mirona1]} and
studied recently by J. Szilasi in \cite{[szilasi]}.

\section{Lagrange spaces}

The variational problem of a Lagrange space determines a canonic
semispray. In this section we prove that for the canonic semispray
of a Lagrange space, there is a unique nonlinear connection that
is metric and it is compatible with the symplectic structure of
the space.

Consider $L^{n}=(M, L)$ a Lagrange space. This means that $L:TM
\longrightarrow \mathbb{R}$ is a regular Lagrangian. In other
words, the (2,0)-type, symmetric, d-tensor field with components
\begin{equation} g_{ij}=\frac{1}{2}\frac{\partial^2 L}{\partial
y^i\partial y^j} \label{lmetric} \end{equation} has rank $n$ on
$TM$. The regularity of the Lagrangian $L$ is also equivalent with
the fact that the Cartan two-form \begin{equation}
\omega=\frac{1}{2}d\left(\frac{\partial L}{\partial
y^i}dx^i\right) \label{omega} \end{equation} is a symplectic
structure on $TM$. The variational problem for the Lagrangian $L$
determines the Euler-Lagrange equations:
\begin{equation}
\frac{\partial L}{\partial x^i} -\frac{d}{dt}\left(\frac{\partial
L}{\partial y^i}\right)=0. \label{eleq}
\end{equation}
Under the assumption of regularity for the Lagrangian $L$, the
system of equations (\ref{eleq}) is equivalent with the following
system of SODE:
\begin{equation} \frac{d^2x^i}{dt^2} +
2G^i\left(x,\frac{dx}{dt}\right)=0. \label{Lsode}\end{equation}
The functions $G^i$ are local coefficients of a semispray $S$ on
$TM$, and they are given by:
\begin{equation}
G^i=\frac{1}{4}g^{ik}\left(\frac{\partial ^{2}L}{\partial y^{k}\partial x^{h}}%
y^{h}-\frac{\partial L}{\partial x^{k}}\right). \label{gi}
\end{equation}
We refer to this semispray as to the canonic semispray of the
Lagrange space. The canonic semispray $S$ of the Lagrange space
$L^n$ is the unique vector field that satisfies the equation
$i_S\omega = -(1/2)d E_L$, where $E_L=y^i({\partial L}/{\partial
y^i})-L$ is the energy of the Lagrange space $L^n$. The semispray
$S$ determines a canonic nonlinear connections, which depends only
on the fundamental function $L$. The local coefficients of this
nonlinear connection are given by \cite{[grifone]}
\begin{equation}
N^i_j=\frac{\partial G^i}{\partial y^j}. \label{nij}
\end{equation}

For the canonic semispray $S$ and the canonic nonlinear connection
$N$ consider $\nabla$ the induced dynamical covariant derivative
(\ref{covariant}).
\begin{theorem} \label{uniquen}
For a Lagrange space $L^n$, the canonic nonlinear connection
(\ref{nij}) is the unique nonlinear connection $N$ that satisfies:
\begin{itemize} \item[1)] The horizontal subbundle $NTM$ is a
Lagrangian subbundle of $TTM$, which means that:
\begin{equation}
\omega(hX, hY)=0, \forall X,Y \in \chi(TM). \label{omegahh}
\end{equation}
\item[2)] The metric tensor $g_{ij}$ of the Lagrange space is
covariant constant with respect to the dynamical covariant
derivative induced by $(S,N)$, which is equivalent to:
\begin{equation}
\nabla g=0. \label{nablag0}
\end{equation}
\end{itemize}
\end{theorem}
\begin{proof} First we prove that conditions
(\ref{omegahh}) and (\ref{nablag0}) uniquely determine a nonlinear
connection. Then, we show that this nonlinear connection is the
canonic nonlinear connection of the Lagrange space.

Consider $N$ a nonlinear connection with local coefficients
$N^i_j$. We want to express the symplectic form $\omega$ using the
adapted cobasis $\{d x^i, \delta y^i\}$. If we use expression
(\ref{omega}) and replace $d y^i=\delta y^i - N^i_jd x^j$ we
obtain:
\begin{equation}
\begin{array}{l}
\omega =g_{ij} (\delta y^j-N^j_k d x^k)\wedge d x^i +
\displaystyle\frac{1}{4}\left(\displaystyle\frac{\partial^2
L}{\partial y^i\partial x^j} - \displaystyle\frac{\partial^2
L}{\partial x^i
\partial y^j}\right)d x^j\wedge d x^i
\vspace{2mm}\\
= g_{ij}\delta y^j\wedge d x^i +
\displaystyle\frac{1}{2}\left[-N_{ij} + N_{ji}+
\displaystyle\frac{1}{2}\left(\displaystyle\frac{\partial^2
L}{\partial y^i\partial x^j} - \displaystyle\frac{\partial^2
L}{\partial x^i
\partial y^j}\right)\right]d x^j\wedge d x^i
\vspace{2mm}\\
= g_{ij}\delta y^j\wedge d x^i + \left[-N_{[ij]} +
\displaystyle\frac{1}{4}\left(\displaystyle\frac{\partial^2
L}{\partial y^i\partial x^j} - \displaystyle\frac{\partial^2
L}{\partial x^i
\partial y^j}\right)\right]d x^j\wedge d x^i,
\end{array}
\label{Lomegaad} \end{equation} where $N_{[ij]}$ denotes the skew
symmetric part of $N_{ij}:=g_{ik}N^k_j$. We have that
(\ref{omegahh}) is true if and only if the second term of the
right hand side of (\ref{Lomegaad}) vanishes. Consequently, we
have that (\ref{omegahh}) is true if and only if:
\begin{equation}
N_{[ij]}=\frac{1}{2}(N_{ij}-N_{ji})=
\displaystyle\frac{1}{4}\left(\displaystyle\frac{\partial^2
L}{\partial y^i\partial x^j} - \displaystyle\frac{\partial^2
L}{\partial x^i
\partial y^j}\right). \label{Lanij}
\end{equation}
Expression (\ref{Lanij}) tells us that the skew symmetric part of
$N_{ij}$ is uniquely determined by condition (\ref{omegahh}) and
hence $N_{[ij]}$ is uniquely determined by the symplectic
structure $\omega$. Next, we prove that the symmetric part of
$N_{ij}$ is perfectly determined by metric condition
(\ref{nablag0}). In local coordinates, condition (\ref{nablag0})
is equivalent to:
\begin{equation}
S(g_{ij}) = g_{mj}N^m_i + g_{im}N^m_j=N_{ij}+N_{ji}=2N_{(ij)}.
\label{Lhelmotz2}
\end{equation}
Expressions (\ref{Lanij}) and (\ref{Lhelmotz2}) uniquely determine
the local coefficients $N^i_j$ of the nonlinear connection $N$
that satisfies (\ref{omegahh}) and (\ref{nablag0}). These
coefficients are given by:
\begin{equation}
\begin{array}{ll}
N^i_j & = g^{ik}N_{kj}=g^{ik}(N_{(kj)}+N_{[kj]}) \vspace{2mm}\\
& = \displaystyle\frac{1}{2}g^{ik}\left[ S(g_{kj}) +
\displaystyle\frac{1}{2}\left(\displaystyle\frac{\partial^2
L}{\partial y^k\partial x^j} - \displaystyle\frac{\partial^2
L}{\partial x^k \partial y^j}\right)\right].
\end{array} \label{uniquenij} \end{equation}
We prove now that the nonlinear connection (\ref{nij}) of a
Lagrange space is the unique one that satisfies (\ref{omegahh})
and (\ref{nablag0}). For this we have to show that the nonlinear
connection with local coefficients (\ref{nij}) satisfies
(\ref{uniquenij}). The coefficients $N^i_j$ of the canonic
nonlinear connection (\ref{nij}) of a Lagrange space can be
written as:
$$
\begin{array}{ll}
N^i_j=\displaystyle\frac{\partial G^i}{\partial y^j} = &
\displaystyle\frac{1}{4}\displaystyle\frac{\partial
g^{ip}}{\partial
y^j}\left(\displaystyle\frac{\partial^2L}{\partial y^p\partial
x^m}y^m
-\displaystyle\frac{\partial L^2}{\partial x^p}\right)  \vspace{2mm}\\
& + \displaystyle\frac{1}{4} g^{ip}
\left(2\displaystyle\frac{\partial g_{jp}}{\partial x^m} y^m -
\displaystyle\frac{\partial^2 L}{\partial y^j\partial x^p}\right)
+ \displaystyle\frac{1}{4}g^{ip}\displaystyle\frac{\partial^2
L}{\partial y^p\partial x^j}.
\end{array}
$$
If we multiply the above formula by $g_{is}$ we obtain:
$$
N_{sj}:=g_{si}N^i_j=-\frac{\partial g_{is}}{\partial y^j}G^i +
\frac{1}{2}\frac{\partial g_{sj}}{\partial x^i} y^i +
\frac{1}{4}\left(\frac{\partial^2 L}{\partial y^s\partial x^j} -
\frac{\partial^2 L}{\partial x^s \partial y^j}\right),
$$
which is equivalent to:
\begin{equation}
N_{ij}=\frac{1}{2}S(g_{ij})+ \frac{1}{4}\left(\frac{\partial^2
L}{\partial y^i\partial x^j} - \frac{\partial^2 L}{\partial x^i
\partial y^j}\right). \label{Lnij}
\end{equation}
We can see that expressions (\ref{Lnij}) and (\ref{uniquenij}) are
equivalent to one another. From expression (\ref{Lnij}) it follows
that the canonic nonlinear connection of a Lagrange space with
local coefficients (\ref{nij}) satisfies the two axioms of the
theorem.
\end{proof}

Theorem \ref{uniquen} shows that the canonic nonlinear connection
of a Lagrange space has:
\begin{itemize}
\item[1)] the skew-symmetric part $N_{[ij]}=(1/2)a_{ij}$ uniquely
determined by the symplectic form $\omega=g_{ij}\delta y^j\wedge
dx^i + (1/2)a_{ij} dx^j\wedge dx^i$. \item[2)] the symmetric part
$N_{(ij)}=S(g_{ij})$ uniquely determined by the semispray $S$ and
the metric tensor $g$.
\end{itemize}

Consequently, we can generalize Theorem \ref{uniquen} as follows:
\begin{theorem} Consider $S$ a semispray and $\omega$ a symplectic
structure on $TM$ for which the vertical subbundle $VTM$ is a
Lagrangian subbundle. There exists a unique nonlinear connection
$N$ on $TM$ such that:
\begin{itemize}
\item[1)] The horizontal subbundle $NTM$ is a Lagrangian subbundle
of $TTM$, which means that
$$
\omega(hX,hY)=0, \forall X,Y \in \chi(TM).
$$
\item[2)] The metric tensor $g_{ij}$ of the generalized Lagrange
space is covariant constant, which means that $$ \nabla g=0.$$
\end{itemize}
\end{theorem}

With respect to the adapted cobasis $\{d x^i, \delta y^i=d y^i +
N^i_j d x^j\}$ of the canonic nonlinear connection the symplectic
form $\omega$ of a Lagrange space $L^n$ has a simple form:
\begin{equation}
\omega = g_{ij} \delta y^j\wedge d x^i. \label{Ladaptedomega}
\end{equation}
Expression (\ref{Ladaptedomega}) is equivalent to (\ref{omegahh}),
which says that symplectic form $\omega$ vanishes if both of its
arguments are horizontal vector fields. One can see also from
expression (\ref{Ladaptedomega}) that $\omega(X,Y)=0$ if both
vectors $X$ and $Y$ are vertical vector fields. Therefore both
horizontal and vertical subbundles are Lagrangian subbundles for
the tangent bundle $TTM$.

A nonlinear connection is perfectly determined by an almost
complex structure $\mathbb{F}$ given by:
\begin{equation}
\mathbb{F}=\frac{\delta}{\delta x^i}\otimes \delta y^i -
\frac{\partial}{\partial y^i}\otimes d x^i. \label{fstructure}
\end{equation}
For a GL-metric $g$ one can use a nonlinear connection $N$ to
define a nondegenerate metric tensor $\mathbb{G}$ on $TM$ that
preserves the horizontal and vertical distributions:
\begin{equation}
\mathbb{G}=g_{ij}dx^i\otimes dx^j + g_{ij}\delta y^i \otimes
\delta y^j. \label{metric}
\end{equation}
The compatibility condition (\ref{omegahh}) implies
$$
\omega(X,Y)=\mathbb{G}(\mathbb{F}X,Y) \textrm{ and } \mathbb{G}(X,
Y)=\mathbb{G}(\mathbb{F}Y, X), \forall X, Y \in \chi(TM).
$$
Consequently, the pair $(\mathbb{G}, \mathbb{F})$ is an almost
Hermitian structure on $TM$.

Theorem \ref{uniquen} shows that the canonic nonlinear connection
of a Lagrange space is metric. Next, as we did for a GL-metric, we
can determine the family of all metric nonlinear connections for a
Lagrange space. For this we do not require anymore the
compatibility of the nonlinear connection with the symplectic
structure.
\begin{proposition}
The family of all nonlinear connections that are compatible with
the metric tensor of a Lagrange space is given by:
\begin{equation}
N^i_j=N^{ci}_j + O^{ki}_{jm} X^m_k. \label{Lmetricnon}
\end{equation}
Here $X^m_k$ is an arbitrary (1,1)-type d-tensor field,
$O^{ki}_{jm}$ is the Obata operator (\ref{GLobata}) and $N^{ci}_j$
are the local coefficients of the canonic nonlinear connection of
the Lagrange space.
\end{proposition}
The proof of this result is similar with that of Theorem
\ref{family}.

\section{Discussions}

In this work we start with a system of second order differential
equations and a generalized Lagrange metric and we determine a
metric nonlinear connection.

This is a different approach of the inverse problem of Lagrangian
Mechanics, where for a given system of SODE and an associated
nonlinear connection, we seek for a metric tensor that makes the
nonlinear connection metric.

The metric nonlinear connection we determine in Theorem 2.2 is not
unique. However its symmetric part is uniquely determined by the
metric compatibility. The metric nonlinear connection given by
expression (10) depends on both structures: semispray and
generalized Lagrange metric. Hence, this nonlinear connection is
different from the nonlinear connection, given by expression (20),
which one usually associate to a semispray. The two nonlinear
connections coincide for the particular case when the metric
tensor is Lagrangian.

A metric nonlinear connection is uniquely determine if we add a
condition that determines its skew-symmetric part. This can be
done if we require the compatibility of the nonlinear connection
with a symplectic structure as we did in Theorem 3.2. For a
Lagrange space we prove in Theorem 3.1 that the metric structure
and the symplectic structure uniquely determine the nonlinear
connection we associate with Euler-Lagrange equations.

\vspace*{3mm}

\noindent\textbf{Acknowledgment} This work has been partially
supported by CNCSIS grant from the Romanian Ministry of Education.

\end{document}